# Solution of Non-Square Fuzzy Linear Systems

Nizami Gasilov [1], Afet Golayoğlu Fatullayev [2] and Şahin Emrah Amrahov [3]

[1] *Baskent University, Ankara, 06810 Turkey; phone: (+90)3122341010/1220; fax: (+90)3122341051*
*E-mail: gasilov@baskent.edu.tr (corresponding author)*
[2] *Baskent University, Ankara, 06810 Turkey; phone: (+90)3122341010/1790; fax: (+90)3122341179*
*E-mail: afet@baskent.edu.tr*
[3] *Ankara University, Computer Engineering Department, Ankara, 06500 Turkey; phone: (+90)3123800026/162*
*E-mail: emrah@eng.ankara.edu.tr*

In this paper, a linear system of equations with crisp coefficients and fuzzy right-hand sides is investigated. All possible cases pertaining to the number of variables, *n*, and the number of equations, *m*, are dealt with. A solution is sought not as a fuzzy number vector, as usual, but as a fuzzy set of vectors. Each vector in the solution set solves the given fuzzy linear system with a certain possibility. Assuming that the coefficient matrix is a full rank matrix, three cases are considered: For $m = n$ (square system), the solution set is shown to be a parallelepiped in coordinate space and is expressed by an explicit formula. For $m > n$ (overdetermined system), the solution set is proved to be a convex polyhedron and a novel geometric method is proposed to compute it. For $m < n$ (underdetermined system), by determining the contribution of free variables, general solution is computed. From the results of three cases mentioned above, a method is proposed to handle the general case, in which the coefficient matrix is not necessarily a full rank matrix. Comprehensive examples are provided and investigated in depth to illustrate each case and suggested method.

## 1. INTRODUCTION

Linear algebraic equation systems are commonly used for modeling application problems. They play important role for solving many problems in science and technology. When a problem includes uncertain data it is more suitable to represent them by fuzzy quantities rather than crisp ones. Therefore we can model such problems with fuzzy linear systems. Recently, there have been significant researches on fuzzy linear systems (FLS) [1-12], which are based mainly on the work of Friedman et al [1]. In these researches, the system is usually considered to be square and consistent. For overdetermined and underdetermined systems there are very few studies such as [13-17]. In these studies Asady et al. [13], Allahviranloo and Kermani [14], and Zheng and Wang [15] have extended an embedding method given by Wu and Ma [18, 19] for solving $m \times n$ fuzzy linear systems. In [17] the same idea is used to solve the fuzzy linear matrix equation $AXB = C$.

In researches, mentioned above, the solution is sought to be a fuzzy number vector. Suggested methods are valid only for coefficient matrix satisfying certain conditions. This is because the existence theorem of Friedman takes place only in the case when the matrix can be represented as the product of a permutation matrix and a diagonal matrix. To fill this gap,

Gasilov et al [20, 21] suggested a new approach for linear systems, which is a geometric approach based on the properties of linear transformations. Unlike the other approaches, in this approach, the solution is considered to be a fuzzy set of vectors. Therefore, the suggested approach can be applied to all linear square systems. The most important difference of the suggested approach from the others will be explained in details later where the concept of the solution will be discussed.

In this paper, we consider an $m \times n$ system with crisp coefficient matrix and with a vector of fuzzy numbers on the right-hand side. We generalize the geometric approach, proposed in [20, 21] for the square systems, for the non-square cases. Therefore, there is no additional requirement on the system matrix, which makes this approach superior to other approaches on the subject. We look for the solution as a fuzzy set. Each vector in the set satisfies the system with a certain possibility. In [22], the solutions of the interval linear systems has been classified as united solution set (USS), tolerable solution set (TSS) and controllable solution set (CSS). In our approach, the solution is in the sense of USS. We first determine the solution set when the coefficient matrix is a full rank matrix. For $n \times n$ square systems we use the geometric approach proposed in [20, 21]. For underdetermined systems, we identify the contribution of free variables to the general solution. For overdetermined systems, we propose a new geometric method. We show that the solution set is a convex polyhedron and we present an algorithm of finding it. In the particular case when the system is square, the convex polyhedron turns into a parallelepiped. For each point of the polyhedron, we determine the possibility that the point satisfies the system. Finally, we suggest a method to handle $m \times n$ fuzzy systems, where the coefficient matrix needs not be a full rank matrix.

This paper consists of 7 sections including the Introduction. Fuzzy linear systems are defined in Section 2. In Section 3, full-rank square systems are considered. Underdetermined systems with full rank matrix are studied in Section 4. In Section 5, a solution method is presented for overdetermined systems with full rank matrix. The general case is considered in Section 6. Conclusions are drawn in Section 7.

## 2. FUZZY LINEAR SYSTEMS

Below, we use the notation $\tilde{f} = (\underline{f}(r), \overline{f}(r))$, $(0 \leq r \leq 1)$ to indicate a fuzzy number in parametric form. We denote $\underline{\underline{f}} = \underline{f}(0)$ and $\overline{\overline{f}} = \overline{f}(0)$ to indicate the left and the right limits of $\tilde{f}$, respectively. We represent a triangular fuzzy number as $\tilde{f} = (l, m, r)$, for which we have $\underline{\underline{f}} = l$ and $\overline{\overline{f}} = r$.

Let $a_{ij}$, $(1 \leq i \leq m, 1 \leq j \leq n)$ be crisp numbers and $\tilde{f}_i = (\underline{f_i}(r), \overline{f_i}(r))$, $(0 \leq r \leq 1, 1 \leq i \leq m)$ be fuzzy numbers. The system

$$\begin{cases} a_{11}x_1 + a_{12}x_2 + \ldots + a_{1n}x_n = \tilde{f}_1 \\ a_{21}x_1 + a_{22}x_2 + \ldots + a_{2n}x_n = \tilde{f}_2 \\ \vdots \\ a_{m1}x_1 + a_{m2}x_2 + \ldots + a_{mn}x_n = \tilde{f}_m \end{cases} \quad (1)$$

is called an $m \times n$ fuzzy linear system.

One can rewrite (1) as follows using matrix notation:

$$A\tilde{X} = \tilde{B} \quad (2)$$

where $A = [a_{ij}]$ is an $m \times n$ crisp matrix and $\tilde{B} = (\tilde{f}_1, \tilde{f}_2, \ldots, \tilde{f}_m)^T$ is a vector of fuzzy numbers.

One can look for solutions of fuzzy linear systems in two ways: solutions as fuzzy number vectors or solutions as fuzzy sets formed by vectors. In the first way, a fuzzy number vector $\tilde{X} = (\tilde{x}_1, \tilde{x}_2, \ldots, \tilde{x}_n)^T$ is called to be a solution of the system if the fuzzy numbers $\tilde{x}_1, \tilde{x}_2, \ldots, \tilde{x}_n$ satisfy each equation of the system (1).

Fuzzy number vectors are convenient tool for expressing the solution in a simple and effective form but it gives rise to some difficulties.

First of all, fuzzy number vector solution exists for only some special matrices (Friedman et al., 1998, 2000).

For example, the system

$$\begin{bmatrix} 1 & 2 \\ 0 & 1 \end{bmatrix} \begin{bmatrix} x \\ y \end{bmatrix} = \begin{bmatrix} (-2, -1, 0) \\ (\ 0,\ 1,\ 2) \end{bmatrix} \qquad (3)$$

has no fuzzy number solution.

To see why, let fuzzy number solution be $\tilde{x} = (\underline{x}, x_{cr}, \overline{x})$ and $\tilde{y} = (\underline{y}, y_{cr}, \overline{y})$. Then from the 2nd equation we have immediately: $\tilde{y} = (0, 1, 2)$. We put this value into the 1st equation and obtain: $(\underline{x}, x_{cr}, \overline{x}) + 2 \cdot (0, 1, 2) = (-2, -1, 0)$. From here we have $\underline{x} = -2$ and $\overline{x} = -4$. But this contradicts with that the right limit of a fuzzy number must be greater or equal to the left one. So, the system (3) has no fuzzy number solution. On the other hand, it can be seen that, for instance, the vector $(x, y) = (-3, 1)$ (the crisp solution) with membership degree $\mu = 1$ solves both equations of the system.

The second difficulty with fuzzy number vector solution is that, even if such a solution exists, it may not expose the solution adequately. The following example explains this situation.

**Example 1.** Solve the system

$$\begin{bmatrix} 1 & -2 \\ 1 & 3 \end{bmatrix} \begin{bmatrix} x \\ y \end{bmatrix} = \begin{bmatrix} (-2, 1, 4) \\ (2, 6, 10) \end{bmatrix} \qquad (4)$$

This system has fuzzy number vector solution $(\tilde{x}, \tilde{y})$ with components $\tilde{x} = (2, 3, 4)$ and $\tilde{y} = (0, 1, 2)$. The vector $(x, y) = (5, 1)$ does not belong to the set expressed by the vector $(\tilde{x}, \tilde{y})$. On the other hand, it can be seen that the vector $(x, y) = (5, 1)$ with membership degree $\mu = \frac{1}{3}$ is a solution of the system:

$1 \cdot 5 - 2 \cdot 1 = 3 \in (-2, 1, 4); \qquad \mu_1 = \mu_{(-2, 1, 4)}(3) = \frac{1}{3};$

$1 \cdot 5 + 3 \cdot 1 = 8 \in (2, 6, 10); \qquad \mu_2 = \mu_{(2, 6, 10)}(8) = \frac{1}{2}$

Therefore we have $\mu = \min\{\mu_1, \mu_2\} = \frac{1}{3}$.

We will give the solution of the above system (4) later, using the method that we proposed (See, Section 3).

Summarizing, fuzzy number vector is not sufficient instrument to express the solution.

To overcome the difficulties caused by fuzzy number vector solutions, in this paper, we seek the solution $\tilde{X}$ as a fuzzy set of $n$-dimensional vectors:

$$\widetilde{X} = \{\mathbf{x} = (x_1, x_2, \ldots, x_n) \mid \underline{f_i} \leq a_{i1}x_1 + a_{i2}x_2 + \ldots + a_{in}x_n \leq \overline{\overline{f_i}}, \quad 1 \leq i \leq n\}$$

In other words, we call a vector $\mathbf{x} = (x_1, x_2, \ldots, x_n)$ to be an element of the fuzzy solution set, if it satisfies each $i$-th equation for some $z_i \in \widetilde{f_i}$. The possibility, or membership, of the vector $\mathbf{x} = (x_1, x_2, \ldots, x_n)$ in the solution set $\widetilde{X}$ is determined by the equation satisfied with least possibility. It means that

$$\mu_{\widetilde{X}}(\mathbf{x}) = \min\{\mu_1(\mathbf{x}), \mu_2(\mathbf{x}), \ldots, \mu_m(\mathbf{x})\}, \tag{5}$$

where $\mu_i(\mathbf{x}) = \mu_{\widetilde{f_i}}(a_{i1}x_1 + a_{i2}x_2 + \ldots + a_{in}x_n)$. That is, $\mu_i(\mathbf{x})$ is the possibility of the value obtained in the left-hand side of $i$-th equation with respect to the fuzzy number $\widetilde{f_i}$.

In regard to matrix $A$, we shall consider the following cases:

1) $m = n$ and $rank(A) = n$ (square system with full rank matrix);
2) $n > m$ and $rank(A) = m$ (underdetermined system with full rank matrix);
3) $m > n$ and $rank(A) = n$ (overdetermined system with full rank matrix);
4) $rank(A) = k \leq \min\{n, m\}$ (general case).

Let $k$ be the rank of $A$, $k = rank(A)$. We shall assume, without loss of generality, that the linearly independent rows and columns of $A$ reside in the upper-left corner of $A$. The $k \times k$ square submatrix in this place we shall denote by $K$. We have $rank(A) = rank(K) = k$ and $K$ is invertible.

Let $\mathbf{y} \in R^k$ and $\mathbf{z} \in R^{n-k}$ be the vectors of leading and free variables, respectively. Then

$$\mathbf{x} = \begin{bmatrix} \mathbf{y} \\ \cdots \\ \mathbf{z} \end{bmatrix}$$

(here and below we use $\cdots$ and $\vdots$ as delimiters inside matrices).

We can represent $A$ by $A = [L \vdots R]$, where $L$ and $R$ are coefficient matrices of leading and free variables, respectively. $L$ is an $m \times k$ and $R$ is an $m \times (n-k)$ matrices. In the cases, where we refer to $K$, we shall represent $A$ as $A = \begin{bmatrix} K & \vdots & G \\ \cdots & \vdots & \cdots \\ M & \vdots & H \end{bmatrix}$. Here $M$ is an $(m-k) \times k$; $G$ is an $k \times (n-k)$ and $H$ is an $(m-k) \times (n-k)$ matrices.

## 3. SQUARE SYSTEM WITH FULL RANK MATRIX

This is the most curious and investigated case, in which $A$ is a square matrix and invertible (i.e. $m = n$ and $rank(A) = n$). Regarding this case, we could summarize results of [20, 21], below. The main approach in these works is as follows: The vector $\widetilde{B}$ of fuzzy numbers forms a rectangular prism in coordinate space. The image of this prism under linear transformation $T(\mathbf{b}) = A^{-1}\mathbf{b}$ is a parallelepiped, which determines the solution set $\widetilde{X}$.

It is convenient to represent the right-hand side vector as $\tilde{B} = \mathbf{b}_{cr} + \tilde{\mathbf{b}}$, where $\mathbf{b}_{cr}$ is a vector with possibility 1 and denotes the vertex (crisp part) of the fuzzy region $\tilde{B}$, while $\tilde{\mathbf{b}}$ denotes the uncertain part the vertex of which is at the origin.

The vector $\mathbf{x}_{cr} = A^{-1}\mathbf{b}_{cr}$ represents the crisp solution.

In the case when right-hand side of (2) consists of triangular fuzzy numbers $\tilde{f}_i = (l_i, m_i, r_i)$, we have $(b_{cr})_i = m_i$ and $\tilde{b}_i = (l_i - m_i, 0, r_i - m_i) = (\underline{\underline{b_i}}, 0, \overline{\overline{b_i}})$. Let $\mathbf{e}_1, \mathbf{e}_2, \ldots, \mathbf{e}_n$ be standard basis vectors. Then the rectangular prism, defined by $\tilde{B}$, can be represented as follows:

$$\tilde{B} = \{\mathbf{b} = \mathbf{b}_{cr} + c_1\mathbf{e}_1 + c_2\mathbf{e}_2 + \ldots + c_n\mathbf{e}_n \mid \underline{\underline{b_i}} \leq c_i \leq \overline{\overline{b_i}}\}.$$

An $\alpha$-cut of $\tilde{B}$ can be expressed in following way:

$$B_\alpha = \{\mathbf{b} = \mathbf{b}_{cr} + c_1\mathbf{e}_1 + c_2\mathbf{e}_2 + \ldots + c_n\mathbf{e}_n \mid (1-\alpha)\underline{\underline{b_i}} \leq c_i \leq (1-\alpha)\overline{\overline{b_i}}\}. \tag{6}$$

Let $\mathbf{g}_i = A^{-1}\mathbf{e}_i$. It can be seen from (6) that an $\alpha$-cut of the solution and the solution itself are represented as follows:

$$X_\alpha = \{\mathbf{x} = \mathbf{x}_{cr} + c_1\mathbf{g}_1 + c_2\mathbf{g}_2 + \ldots + c_n\mathbf{g}_n \mid (1-\alpha)\underline{\underline{b_i}} \leq c_i \leq (1-\alpha)\overline{\overline{b_i}}\}; \tag{7}$$

$$\tilde{X} = X_0 \text{ with } \mu_{\tilde{X}}(\mathbf{x}) = 1 - \max_{1\leq i\leq n} \gamma_i, \text{ where } \gamma_i = \begin{cases} c_i/\overline{\overline{b_i}}, & c_i \geq 0 \\ c_i/\underline{\underline{b_i}}, & c_i < 0 \end{cases} \tag{8}$$

In the general case when right-hand side consists of parametric fuzzy numbers $\tilde{f}_i = (\underline{f_i}(r), \overline{f_i}(r))$ the solution can be represented as follows:

$$X_\alpha = \{\mathbf{x} = \mathbf{x}_{cr} + c_1\mathbf{g}_1 + c_2\mathbf{g}_2 + \ldots + c_n\mathbf{g}_n \mid \underline{f_i}(\alpha) - (b_{cr})_i \leq c_i \leq \overline{f_i}(\alpha) - (b_{cr})_i\}; \tag{9}$$

$$\tilde{X} = X_0 \text{ with } \mu_{\tilde{X}}(\mathbf{x}) = \min_{1\leq i\leq n} \alpha_i, \tag{10}$$

where

$$\alpha_i = \begin{cases} \overline{f_i}^{-1}(k_i), & k_i > \overline{f_i}(1) \\ 1, & \underline{f_i}(1) \leq k_i \leq \overline{f_i}(1) \text{; and } k_i = (b_{cr})_i + c_i \\ \underline{f_i}^{-1}(k_i), & k_i < \underline{f_i}(1) \end{cases} \tag{11}$$

Note that $\mathbf{g}_i = A^{-1}\mathbf{e}_i$ is determined by the $i$-th column of matrix $A^{-1}$, so no additional calculations required to construct the parallelepipeds $X_\alpha$ and $\tilde{X} = X_0$.

To demonstrate the solution method we consider Example 1 again.
**Example 1.** Solve the system
$$\begin{bmatrix} 1 & -2 \\ 1 & 3 \end{bmatrix}\begin{bmatrix} x \\ y \end{bmatrix} = \begin{bmatrix} (-2, 1, 4) \\ (2, 6, 10) \end{bmatrix}$$
**Solution.**

We have $A^{-1} = \frac{1}{5}\begin{bmatrix} 3 & 2 \\ -1 & 1 \end{bmatrix}$

The vector $\tilde{B} = (\tilde{a}, \tilde{b})$, the components of which are triangular fuzzy numbers $\tilde{a} = (-2, 1, 4)$ and $\tilde{b} = (2, 6, 10)$ on the right-hand side of the system, forms a rectangular region in the coordinate plane. This region is shown by the rectangle $ABCD$ in Fig. 1. The image of this rectangle under linear transformation $T(\mathbf{b}) = A^{-1}\mathbf{b}$ is the parallelogram $A'B'C'D'$, which determines the solution set $\tilde{X}$.

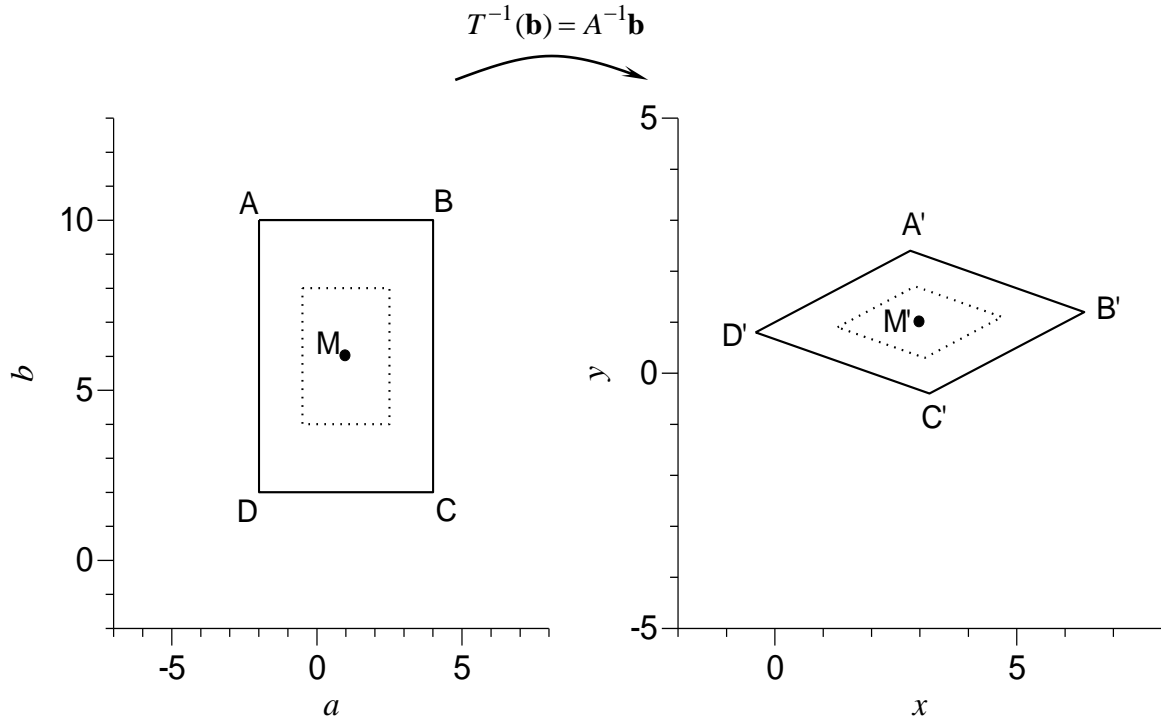

**Figure 1**. Fuzzy numbers $\tilde{a}$ and $\tilde{b}$ on the right-hand side of the system form a rectangular region $ABCD$. The linear transformation $T(\mathbf{b}) = A^{-1}\mathbf{b}$ maps this rectangle to the parallelogram $A'B'C'D'$, which determines the solution set $\tilde{X}$.

Using the formulas (7)-(8), we can express the solution $\tilde{X}$ formally:

$$\tilde{X} = \left\{ \begin{bmatrix} x \\ y \end{bmatrix} = \begin{bmatrix} 3 \\ 1 \end{bmatrix} + c_1 \begin{bmatrix} 3 \\ -1 \end{bmatrix} + c_2 \begin{bmatrix} 2 \\ 1 \end{bmatrix} \;\middle|\; \begin{array}{l} c_1 \in [-3, 3], \\ c_2 \in [-4, 4] \end{array} \right\}$$

*Remark*: One can easily check that the vector $(x, y) = (5, 1)$ mentioned above (See, Example 1 in Section 2) belongs to the set $\tilde{X}$ obtained by the proposed method (it corresponds to the values $c_1 = 0.4$, $c_2 = 0.4$ of the coefficients). This example can be considered as a justification of the fact that the proposed method is more appropriate in general.

## 4. UNDERDETERMINED SYSTEM WITH FULL RANK MATRIX

In this section we assume that $n > m$ and $rank(A) = m$.
We first consider a crisp system and then adopt the results for a fuzzy system. We should note that a system has infinitely many solutions in the considered case.

### 4.1. Solution of crisp underdetermined systems

Consider the crisp system $A\mathbf{x} = \mathbf{b}$. We have $k = rank(A) = m$, $\mathbf{x} = \begin{bmatrix} \mathbf{y} \\ \cdots \\ \mathbf{z} \end{bmatrix}$ and $A = [K \vdots G]$.

We shall determine the general solution of the system:
$$A\mathbf{x} = \mathbf{b} \Rightarrow K\mathbf{y} + G\mathbf{z} = \mathbf{b} \Rightarrow \mathbf{y} = K^{-1}\mathbf{b} - K^{-1}G\mathbf{z}$$

Taking $\mathbf{z} = \mathbf{p}$, where $\mathbf{p} \in R^{n-m}$ is a vector of parameters (arbitrary values such as $s, t, \ldots$, assigning to the free variables), we could write the general solution as:

$$\mathbf{x} = \begin{bmatrix} \mathbf{y} \\ \cdots \\ \mathbf{z} \end{bmatrix} = \begin{bmatrix} K^{-1}\mathbf{b} - K^{-1}G\mathbf{p} \\ \cdots \\ \mathbf{p} \end{bmatrix}$$

or

$$\mathbf{x} = \begin{bmatrix} \mathbf{y} \\ \cdots \\ \mathbf{z} \end{bmatrix} = \begin{bmatrix} K^{-1}\mathbf{b} \\ \cdots \\ \mathbf{0} \end{bmatrix} + \begin{bmatrix} -K^{-1}G \\ \cdots \\ I \end{bmatrix} \mathbf{p}, \quad (12)$$

where $I$ is the $(n-m) \times (n-m)$ identity matrix.

In the last formula, the solution is presented as sum of the solution of non-homogeneous system with $\mathbf{z} = \mathbf{0}$ and the general solution of associated homogeneous system.

**Example 2.** Solve $\begin{cases} -x + 2y + 3z = 1 \\ 3x + 4y - 2z = 17 \end{cases}$

**Solution.** We establish the values of parameters indicated above and compute the solution by the formula (12).

$$A = [K \vdots G] = \begin{bmatrix} -1 & 2 & \vdots & 3 \\ 3 & 4 & \vdots & -2 \end{bmatrix}; \quad \mathbf{b} = \begin{bmatrix} 1 \\ 17 \end{bmatrix}; \quad K = \begin{bmatrix} -1 & 2 \\ 3 & 4 \end{bmatrix}$$

Therefore,
$$K^{-1} = \tfrac{1}{10}\begin{bmatrix} -4 & 2 \\ 3 & 1 \end{bmatrix}$$

Hence,
$$K^{-1}\mathbf{b} = \begin{bmatrix} 3 \\ 2 \end{bmatrix}.$$

$$K^{-1}G = \tfrac{1}{10}\begin{bmatrix} -4 & 2 \\ 3 & 1 \end{bmatrix}\begin{bmatrix} 3 \\ -2 \end{bmatrix} = \tfrac{1}{10}\begin{bmatrix} -16 \\ 7 \end{bmatrix}$$

$$-K^{-1}G\mathbf{p} = -\tfrac{1}{10}\begin{bmatrix} -16 \\ 7 \end{bmatrix}[s] = \begin{bmatrix} 1.6s \\ -0.7s \end{bmatrix}$$

General solution: $\begin{bmatrix} x \\ y \\ z \end{bmatrix} = \begin{bmatrix} 3 \\ 2 \\ 0 \end{bmatrix} + s\begin{bmatrix} 1.6 \\ -0.7 \\ 1 \end{bmatrix}; \quad s \in R$

**Example 3.** Solve $\begin{cases} -x+2y-7z+5w=1 \\ 3x+4y+z-5w=17 \end{cases}$

**Solution.**

$$A = [K \vdots G] = \begin{bmatrix} -1 & 2 & \vdots & -7 & 5 \\ 3 & 4 & \vdots & 1 & -5 \end{bmatrix}; \quad \mathbf{b} = \begin{bmatrix} 1 \\ 17 \end{bmatrix}; \quad K = \begin{bmatrix} -1 & 2 \\ 3 & 4 \end{bmatrix}$$

Therefore,

$$K^{-1} = \tfrac{1}{10}\begin{bmatrix} -4 & 2 \\ 3 & 1 \end{bmatrix}$$

Hence,

$$K^{-1}\mathbf{b} = \begin{bmatrix} 3 \\ 2 \end{bmatrix}$$

$$K^{-1}G = \tfrac{1}{10}\begin{bmatrix} -4 & 2 \\ 3 & 1 \end{bmatrix}\begin{bmatrix} -7 & 5 \\ 1 & -5 \end{bmatrix} = \tfrac{1}{10}\begin{bmatrix} 30 & -30 \\ -20 & 10 \end{bmatrix} = \begin{bmatrix} 3 & -3 \\ -2 & 1 \end{bmatrix}$$

$$-K^{-1}G\mathbf{p} = -\begin{bmatrix} 3 & -3 \\ -2 & 1 \end{bmatrix}\begin{bmatrix} s \\ t \end{bmatrix} = \begin{bmatrix} -3s+3t \\ 2s-t \end{bmatrix}$$

General solution: $\begin{bmatrix} x \\ y \\ z \\ w \end{bmatrix} = \begin{bmatrix} 3 \\ 2 \\ 0 \\ 0 \end{bmatrix} + s\begin{bmatrix} -3 \\ 2 \\ 1 \\ 0 \end{bmatrix} + t\begin{bmatrix} 3 \\ -1 \\ 0 \\ 1 \end{bmatrix}; \quad s \in R, \; t \in R$

### 4.2. Solution of fuzzy underdetermined systems

It can be seen that, like to the crisp case, the solution can be presented as sum of the solution of non-homogeneous system with values of free variables taken be equal to zero ($\mathbf{z} = \mathbf{0}$) and the general solution of associated homogeneous system. Note that the first system is fuzzy square system with full rank matrix (which can be solved by the method presented in the previous section), while the second one is a crisp system. Hence, the only difference from the crisp case is that instead of crisp system $K\mathbf{y} = \mathbf{b}$ the fuzzy system $K\mathbf{y} = \tilde{B}$ will be solved. Therefore, the following method can be offered.

*The solution method*: The fuzzy solution set consists of vectors such as $\mathbf{x} = \begin{bmatrix} \mathbf{y} \\ \cdots \\ \mathbf{z} \end{bmatrix}$. Take $\mathbf{z} = \mathbf{p}$, ($\mathbf{p} \in R^{n-m}$). The parallelepiped corresponding to the solution of fuzzy square system $K\mathbf{y} = \tilde{B}$ is determined by (9)-(11) (or (7)-(8) in the case of triangular numbers). This parallelepiped, translated by the vectors $-K^{-1}G\mathbf{z} = -K^{-1}G\mathbf{p}$, gives the $\mathbf{y}$ part of the general solution of the system $A\tilde{X} = \tilde{B}$.

**Example 4.** Solve $\begin{cases} -x+2y+3z = (-1,\,1,\,3) \\ 3x+4y-2z = (15,17,20) \end{cases}$

**Solution.**
We note that the given system is same as in Example 2, but with fuzzy right-hand side.

The solution of the system $K\mathbf{y} = \tilde{B}$, or explicitly,
$$\begin{cases} -x + 2y = (-1, 1, 3) \\ 3x + 4y = (15, 17, 20) \end{cases}$$
is computed using (7)-(8) as follows:
$$X_\alpha = \left\{ \begin{bmatrix} x \\ y \end{bmatrix} = \begin{bmatrix} 3 \\ 2 \end{bmatrix} + c_1 \begin{bmatrix} -0.4 \\ 0.3 \end{bmatrix} + c_2 \begin{bmatrix} 0.2 \\ 0.1 \end{bmatrix} \middle| \begin{array}{l} c_1 \in [-2(1-\alpha), 2(1-\alpha)]; \\ c_2 \in [-2(1-\alpha), 3(1-\alpha)] \end{array} \right\}$$

Also we have:
$$-K^{-1}G\mathbf{p} = -\frac{1}{10}\begin{bmatrix} -16 \\ 7 \end{bmatrix}[s] = \begin{bmatrix} 1.6s \\ -0.7s \end{bmatrix}$$

Based on the computations above, one can determine an $\alpha$-cut of the fuzzy solution set as:
$$X_\alpha = \left\{ \begin{bmatrix} x \\ y \\ z \end{bmatrix} = \begin{bmatrix} 3 \\ 2 \\ 0 \end{bmatrix} + c_1 \begin{bmatrix} -0.4 \\ 0.3 \\ 0 \end{bmatrix} + c_2 \begin{bmatrix} 0.2 \\ 0.1 \\ 0 \end{bmatrix} + s \begin{bmatrix} 1.6 \\ -0.7 \\ 1 \end{bmatrix} \middle| \begin{array}{l} c_1 \in [-2(1-\alpha), 2(1-\alpha)], \\ c_2 \in [-2(1-\alpha), 3(1-\alpha)], \\ s \in R \end{array} \right\}$$

## 5. OVERDETERMINED SYSTEM WITH FULL RANK MATRIX

Now we assume that $m > n$ and $rank(A) = n$.

We note that in the crisp case, if $rank(A) \neq rank(A^\#)$, where $A^\# = [A \vdots \mathbf{b}]$ is the augmented matrix of the system, there is no solution; if $rank(A) = rank(A^\#)$ there exists at least one solution; if $rank(A) = rank(A^\#) = n$ there exists a unique solution.

We propose a new method to solve an overdetermined system with full rank matrix.

*The solution method*:

We first consider the case $n = 2$. Then each equation of the system is in the form $ax + by = \tilde{f}$, where $\tilde{f} = (\underline{f}(r), \overline{f}(r))$ is a fuzzy number, and represents a band bounded by the parallel straight lines $ax + by = \underline{\underline{f}}$ and $ax + by = \overline{\overline{f}}$. Here $\underline{\underline{f}} = \underline{f}(0)$ and $\overline{\overline{f}} = \overline{f}(0)$ are upper and lower values of $\tilde{f}$, respectively. The intersection of all the bands corresponding to the equations of the system gives the solution set $X$. This set is a convex polygon. The vertices of the polygon are some intersection points of the straight lines. But an intersection may not be a vertex. More than two straight lines may meet at a vertex. Based on the foregoing discussion we can apply the following algorithm to determine the vertices and hence the solution set: $m$ equations determine $m$ pairs of straight lines. Straight lines of a particular pair do not cross, as they are parallel. Hence, it is possible only for straight lines belonging to distinct pairs to intersect. The number of combinations of two straight lines taken from distinct pairs are $2^2 C_m^2$. Considering each combination separately, we can determine all intersection points. To decide if an intersection point $(x_*, y_*)$ is a vertex, we go through the list of equations to check if $(x_*, y_*)$ satisfies all of them. If it does, then $(x_*, y_*)$ is a vertex. The convex polygon determined by the vertices is the solution set of the system. We note that a convex $r$-gon with vertices $\mathbf{x}_1, \mathbf{x}_2, \ldots, \mathbf{x}_r$ can be completely specified by:

$$P = \{\mathbf{x}_c + \alpha_1(\mathbf{x}_1 - \mathbf{x}_c) + \alpha_2(\mathbf{x}_2 - \mathbf{x}_c) + \ldots + \alpha_r(\mathbf{x}_r - \mathbf{x}_c) \mid \alpha_i \geq 0; \; \sum_{i=1}^{r} \alpha_i \leq 1\},$$

where $\mathbf{x}_c$ is an interior point of the polygon (for instance, $\mathbf{x}_c = (\mathbf{x}_1 + \mathbf{x}_2 + \ldots + \mathbf{x}_r)/r$).

For $n > 2$, straight line and polygon turn into plane (hyper-plane) and polyhedron, respectively. The vertices are obtained from the intersection of $n$ planes. The number of ways of choosing $n$ planes from $m$ pairs of planes (choosing at most one plane from each pair) is $2^n C_m^n$.

We note that the method described above is valid for the first case ($m = n$ and $rank(A) = n$) as well.

**Example 5.** Solve $\begin{cases} -x + 2y = (-1,\,1,\,3) \\ 3x + 4y = (15, 17, 20) \\ 2x - y = (2,\,3,\,6) \end{cases}$. Check if $\begin{bmatrix} 3 \\ 2 \end{bmatrix}$, $\begin{bmatrix} 2.5 \\ 2 \end{bmatrix}$ and $\begin{bmatrix} 2 \\ 2 \end{bmatrix}$ are solutions and determine their possibility if they are.

**Solution.**

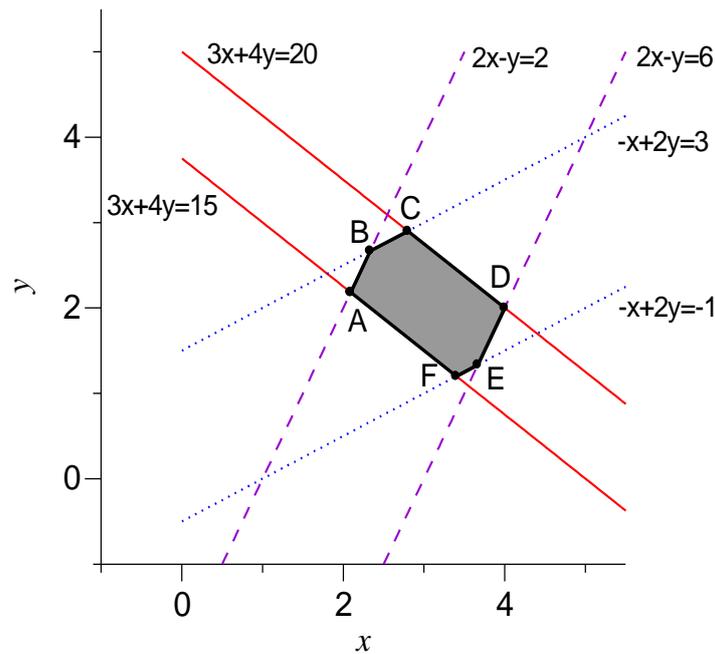

**Figure 2**. Each pair of parallel straight lines bounds a band, determined by an equation of the system. The intersection of all bands gives the solution set $\tilde{X}$ (the polygon *ABCDEF*).

Each equation of the system determines a band in the coordinate plane. The parallel lines which are the boundaries of these bands are shown in Fig. 2. The hexagon formed by the intersection of the bands determine the solution set $\tilde{X}$. The vectors corresponding to the vertices of the hexagon are given by: $\mathbf{x}_A = \mathbf{OA} = (\frac{23}{11}, \frac{24}{11})$, $\mathbf{x}_B = (\frac{7}{3}, \frac{8}{3})$, $\mathbf{x}_C = (2.8, 2.9)$, $\mathbf{x}_D = (4, 2)$, $\mathbf{x}_E = (\frac{11}{3}, \frac{4}{3})$ and $\mathbf{x}_F = (3.4, 1.2)$. We may select an interior point as $\mathbf{x}_c = (\mathbf{x}_A + \mathbf{x}_B + \ldots + \mathbf{x}_F)/6 = (\frac{503}{165}, \frac{1351}{660}) \approx (3.05, 2.05)$. To simplify expression, we shall take $\mathbf{x}_c = (\mathbf{x}_B + \mathbf{x}_E)/2 = (3, 2)$. Then the hexagon is given by:

$$\tilde{X} = \{\mathbf{x}_c + \alpha_1(\mathbf{x}_A - \mathbf{x}_c) + \alpha_2(\mathbf{x}_B - \mathbf{x}_c) + \ldots + \alpha_6(\mathbf{x}_F - \mathbf{x}_c) \mid \alpha_i \geq 0; \sum_{i=1}^{6} \alpha_i \leq 1\},$$

or more explicitly,

$$\tilde{X} = \left\{ \begin{bmatrix} x \\ y \end{bmatrix} = \begin{bmatrix} 3 - \frac{10}{11}\alpha_1 - \frac{2}{3}\alpha_2 - 0.2\alpha_3 + 1\alpha_4 + \frac{2}{3}\alpha_5 + 0.4\alpha_6 \\ 2 + \frac{2}{11}\alpha_1 + \frac{2}{3}\alpha_2 + 0.9\alpha_3 + 0\alpha_4 - \frac{2}{3}\alpha_5 - 0.8\alpha_6 \end{bmatrix} \middle| \begin{array}{l} \alpha_i \geq 0; \\ \sum_{i=1}^{6} \alpha_i \leq 1 \end{array} \right\} \quad (13)$$

The membership function of the fuzzy solution $\tilde{X}$ is displayed in Fig. 3.

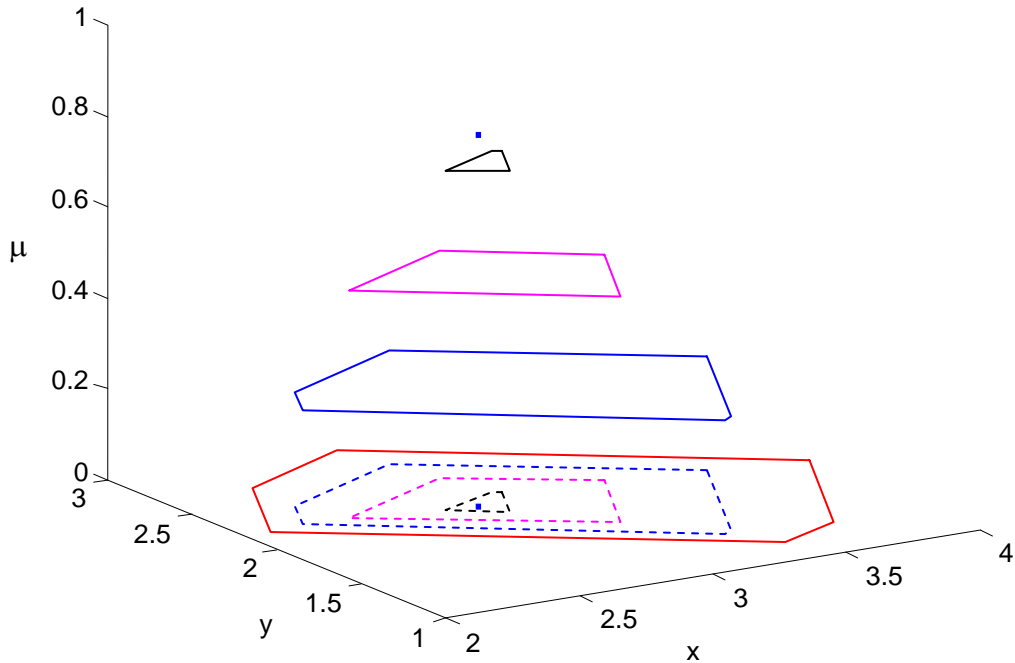

**Figure 3**. The membership function of the fuzzy solution $\tilde{X}$ is represented via $\alpha$-cuts: (from bottom to up) $\alpha \equiv \mu = 0.00, 0.25, 0.50, 0.75$. The solution $\tilde{x} \approx 2.79515$, $\tilde{y} \approx 2.07030$ with highest membership $\alpha \equiv \mu \approx 0.825$ is denoted by a dot.

We shall compute the possibilities of vectors, according to the formula (5):

$\mu_{\tilde{X}}\left(\begin{bmatrix} 3 \\ 2 \end{bmatrix}\right) = \min\{\mu_{\tilde{f}_1}(1), \mu_{\tilde{f}_2}(17), \mu_{\tilde{f}_3}(4)\} = \min\{1, 1, 1 - \frac{4-3}{6-3}\} = \frac{2}{3}$

$\mu_{\tilde{X}}\left(\begin{bmatrix} 2.5 \\ 2 \end{bmatrix}\right) = \min\{\mu_{\tilde{f}_1}(1.5), \mu_{\tilde{f}_2}(15.5), \mu_{\tilde{f}_3}(3)\} = \min\{1 - \frac{1.5-1}{3-1}, \frac{15.5-15}{17-15}, 1\} = 0.25$

$\mu_{\tilde{X}}\left(\begin{bmatrix} 2 \\ 2 \end{bmatrix}\right) = \min\{\mu_{\tilde{f}_1}(2), \mu_{\tilde{f}_2}(14), \mu_{\tilde{f}_3}(2)\} = \mu_{\tilde{f}_2}(14) = 0$. Therefore, $\begin{bmatrix} 2 \\ 2 \end{bmatrix}$ is not a solution.

We note that, for example, if the third equation were $2x - y = (-2, -1, 0)$, the intersection of bands would be empty and, consequently, the system would be inconsistent.

**Example 6.** Solve $\begin{cases} -x + 2y = (-1 + r,\ 3 - 2r^2) \\ 3x + 4y = (15 + r^2,\ 20 - 2\sqrt{r}) \\ 2x - y = (2 + r^3,\ 6 - 3r^2) \end{cases}$. Check if $\begin{bmatrix} 3 \\ 2 \end{bmatrix}$, $\begin{bmatrix} 2.5 \\ 2 \end{bmatrix}$ and $\begin{bmatrix} 2 \\ 2 \end{bmatrix}$ are solutions and determine their possibility if they are.

**Solution.**

We note that the given system is same as in Example 5, but with different right-hand side values, which are represented in parametric form.

For $r = 0$ we have $\begin{bmatrix} (\underline{f_1}(0),\ \overline{f_1}(0)) \\ (\underline{f_2}(0),\ \overline{f_2}(0)) \\ (\underline{f_3}(0),\ \overline{f_3}(0)) \end{bmatrix} = \begin{bmatrix} (-1,\ 3) \\ (15,\ 20) \\ (2,\ 6) \end{bmatrix}$. From here it can be seen that the lower and upper boundary values of right-hand fuzzy numbers are the same as in Example 5. Hence, the solution set $\tilde{X}$ also is the same (Fig. 2), but with different membership function.

We determine the membership functions of the fuzzy numbers on the right-hand side of the system.

$\mu_{\tilde{f}_1}(x) = \begin{cases} x + 1, & -1 \leq x < 0 \\ 1, & 0 \leq x \leq 1 \\ \sqrt{(3-x)/2}, & 1 < x \leq 3 \end{cases}$; $\mu_{\tilde{f}_2}(x) = \begin{cases} \sqrt{x - 15}, & 15 \leq x < 16 \\ 1, & 16 \leq x \leq 18 \\ ((20-x)/2)^2, & 18 < x \leq 20 \end{cases}$;

$\mu_{\tilde{f}_3}(x) = \begin{cases} \sqrt[3]{x - 2}, & 2 \leq x \leq 3 \\ \sqrt{(6-x)/3}, & 3 \leq x \leq 6 \end{cases}$

We can compute the possibilities of given vectors as follows:

$\mu_{\tilde{X}}\left(\begin{bmatrix} 3 \\ 2 \end{bmatrix}\right) = \min\{\mu_{\tilde{f}_1}(1), \mu_{\tilde{f}_2}(17), \mu_{\tilde{f}_3}(4)\} = \min\{1, 1, \sqrt{2/3}\} = \sqrt{2/3} \approx 0.8165$

$\mu_{\tilde{X}}\left(\begin{bmatrix} 2.5 \\ 2 \end{bmatrix}\right) = \min\{\mu_{\tilde{f}_1}(1.5), \mu_{\tilde{f}_2}(15.5), \mu_{\tilde{f}_3}(3)\} = \min\{\sqrt{0.75}, \sqrt{0.5}, 1\} = \sqrt{0.5} \approx 0.7071$

$\mu_{\tilde{X}}\left(\begin{bmatrix} 2 \\ 2 \end{bmatrix}\right) = \min\{\mu_{\tilde{f}_1}(2), \mu_{\tilde{f}_2}(14), \mu_{\tilde{f}_3}(2)\} = \mu_{\tilde{f}_2}(14) = 0$. Therefore $\begin{bmatrix} 2 \\ 2 \end{bmatrix}$ is not a solution.

## 6. GENERAL CASE

Now we consider general case that is $rank(A) = k \leq \min\{n, m\}$.

The general solution of $A\mathbf{x} = \begin{bmatrix} \vdots \\ L & \vdots & R \\ \vdots \end{bmatrix} \begin{bmatrix} \mathbf{y} \\ \cdots \\ \mathbf{z} \end{bmatrix} = \begin{bmatrix} K & \vdots & G \\ \cdots & \vdots & \cdots \\ M & \vdots & H \end{bmatrix} \begin{bmatrix} \mathbf{y} \\ \cdots \\ \mathbf{z} \end{bmatrix} = \tilde{B}$ can be represented as

$\tilde{X} = \tilde{X}_p + X_h$.

Here $\tilde{X}_p$ is the fuzzy solution of the fuzzy non-homogeneous system, with $\mathbf{z} = \mathbf{0}$. We present how one can compute it.

$$A\mathbf{x} = \begin{bmatrix} \vdots & \vdots \\ L & \vdots & R \\ \vdots & \vdots \end{bmatrix} \begin{bmatrix} \mathbf{y} \\ \cdots \\ \mathbf{z} \end{bmatrix} = \tilde{B} \Rightarrow L\mathbf{y} + R\mathbf{z} = \tilde{B} \stackrel{\mathbf{z}=\mathbf{0}}{\Rightarrow} L\mathbf{y} = \tilde{B}.$$

The system $L\mathbf{y} = \tilde{B}$ corresponds to the third case (overdetermined system) if use $k$ instead $n$. One could determine $\mathbf{y}$ by the proposed method above. From this solution, one could construct the set $\tilde{X}_p$ taking $\mathbf{z} = \mathbf{0}$ into account.

$X_h$ is the general solution of the associated crisp homogeneous system:

$$A\mathbf{x} = \begin{bmatrix} K & \vdots & G \\ \cdots & \vdots & \cdots \\ M & \vdots & H \end{bmatrix} \begin{bmatrix} \mathbf{y} \\ \cdots \\ \mathbf{z} \end{bmatrix} = \mathbf{0}.$$

Since $rank(A) = k$, equations in the lower part of the system are linear combinations of the equations in the upper part, hence they could be discarded. Then,

$K\mathbf{y} + G\mathbf{z} = \mathbf{0}$

Therefore,

$\mathbf{y} = -K^{-1}G\mathbf{z}$

Hence $\mathbf{y}_h = -K^{-1}G\mathbf{p}$; $\mathbf{z}_h = \mathbf{p}$. Then,

$$X_h = \left\{ \begin{bmatrix} \mathbf{y}_h \\ \mathbf{z}_h \end{bmatrix} = \begin{bmatrix} -K^{-1}G\mathbf{p} \\ \mathbf{p} \end{bmatrix} \mid \mathbf{p} \in R^{n-k} \right\}$$

**Example 7.** Solve $\begin{cases} -x + 2y - 7z + 5w = (-1,\ 1,\ 3) \\ 3x + 4y\ + z - 5w = (15, 17, 20) \\ 2x\ - y\ + 8z - 7w = (2,\ 3,\ 6) \end{cases}$

Check if $(9,\ 0.5,\ 0,\ 2)^T$, $(3,\ 3,\ 1,\ 1)^T$ and $(-3,\ 4.5,\ 1,\ -1)^T$ are solutions and determine their possibility if they are.

**Solution.**

We have $rank(A) = 2$ and

$$A = \begin{bmatrix} \vdots & \vdots \\ L & \vdots & R \\ \vdots & \vdots \end{bmatrix} = \begin{bmatrix} K & \vdots & G \\ \cdots & \vdots & \cdots \\ M & \vdots & H \end{bmatrix} = \begin{bmatrix} -1 & 2 & \vdots & -7 & 5 \\ 3 & 4 & \vdots & 1 & -5 \\ \cdots & \cdots & \vdots & \cdots & \cdots \\ 2 & -1 & \vdots & 8 & -7 \end{bmatrix}$$

First we determine $\tilde{X}_p$.

$$A\mathbf{x} = \begin{bmatrix} \vdots & \vdots \\ L & \vdots & R \\ \vdots & \vdots \end{bmatrix} \begin{bmatrix} \mathbf{y} \\ \cdots \\ \mathbf{z} \end{bmatrix} = \tilde{B} \Rightarrow L\mathbf{y} + R\mathbf{z} = \tilde{B} \stackrel{\mathbf{z}=\mathbf{0}}{\Rightarrow} L\mathbf{y} = \tilde{B}.$$

The system $L\mathbf{y} = \tilde{B}$ for the current example is the same as that of Example 5. Therefore, we do not have to do additional computation. All we need is to adapt the solutions to 4-dimensional space. Hence the first two components ($x$ and $y$) of the solutions are same as those in Example 5 (see, formula (13)), and the last two components ($z$ and $w$) are 0 (because $\mathbf{z} = \mathbf{0}$).

We shall find the solution of the corresponding homogeneous system:

$$-K^{-1}G = K^{-1}(-G) = \tfrac{1}{10}\begin{bmatrix}-4 & 2\\ 3 & 1\end{bmatrix}\begin{bmatrix}7 & -5\\ -1 & 5\end{bmatrix} = \tfrac{1}{10}\begin{bmatrix}-30 & 30\\ 20 & -10\end{bmatrix} = \begin{bmatrix}-3 & 3\\ 2 & -1\end{bmatrix}$$

$$\begin{bmatrix}x_h\\ y_h\end{bmatrix} = -K^{-1}G\mathbf{p} = \begin{bmatrix}-3 & 3\\ 2 & -1\end{bmatrix}\begin{bmatrix}s\\ t\end{bmatrix} \Rightarrow \begin{bmatrix}x_h\\ y_h\end{bmatrix} = \begin{bmatrix}-3s+3t\\ 2s-t\end{bmatrix}$$

$$X_h = \left\{\begin{bmatrix}x_h\\ y_h\\ z_h\\ w_h\end{bmatrix} = \begin{bmatrix}-3s+3t\\ 2s-t\\ s\\ t\end{bmatrix}\,\middle|\, \begin{array}{l}s\in R;\\ t\in R\end{array}\right\}$$

The general solution $\tilde{X} = \tilde{X}_p + X_h$ is:

$$\tilde{X} = \left\{\begin{bmatrix}x\\ y\\ z\\ w\end{bmatrix} = \begin{bmatrix}3-\tfrac{10}{11}\alpha_1 - \tfrac{2}{3}\alpha_2 - 0.2\alpha_3 + 1\alpha_4 + \tfrac{2}{3}\alpha_5 + 0.4\alpha_6 - 3s + 3t\\ 2+\tfrac{2}{11}\alpha_1 + \tfrac{2}{3}\alpha_2 + 0.9\alpha_3 + 0\alpha_4 - \tfrac{2}{3}\alpha_5 - 0.8\alpha_6 + 2s - t\\ s\\ t\end{bmatrix}\,\middle|\, \begin{array}{l}\alpha_i \geq 0;\ \sum_{i=1}^{6}\alpha_i \leq 1;\\ s\in R;\\ t\in R\end{array}\right\}$$

We shall compute the possibilities of given vectors, according to (5):

$\mu_{\tilde{X}}((9,\ 0.5,\ 0,\ 2)^T) = \min\{\mu_{\tilde{f}_1}(2), \mu_{\tilde{f}_2}(19), \mu_{\tilde{f}_3}(3.5)\} = \min\{0.5, 1-\tfrac{19-17}{20-17}, 1-\tfrac{3.5-3}{6-3}\} = \tfrac{1}{3}$

$\mu_{\tilde{X}}((3,\ 3,\ 1,\ 1)^T) = \min\{\mu_{\tilde{f}_1}(1), \mu_{\tilde{f}_2}(17), \mu_{\tilde{f}_3}(4)\} = \min\{1, 1, 1-\tfrac{4-3}{6-3}\} = \tfrac{2}{3}$

$\mu_{\tilde{X}}((-3,\ 4.5,\ 1,\ -1)^T) = \min\{\mu_{\tilde{f}_1}(0), \mu_{\tilde{f}_2}(15), \mu_{\tilde{f}_3}(4.5)\} = \mu_{\tilde{f}_2}(4.5) = 0 \Rightarrow (-3,\ 4.5,\ 1,\ -1)^T$ is not a solution.

## 7. CONCLUSION

To summarize, in this paper, we investigated a non-square system of $m$ fuzzy linear equations with $n$ variables. The coefficients of the equations are crisp and the right hand sides are fuzzy numbers. Instead of looking for solution to be a vector of fuzzy numbers, we determined fuzzy solution set, consisting of vectors each of that satisfies the system with some possibility. We proposed a new geometric approach to solve an overdetermined system (the case $m > n$). For an underdetermined system (the case $n > m$) we determined the contribution of free variables to general solution. Finally, we suggested a method to compute the fuzzy solution set of the system in general case, which assembles the results of the special cases, mentioned above. Some numerical examples are presented to emphasize different aspects of the suggested method.

For future research we will attempt to apply the proposed method for partially and fully fuzzy linear systems. The geometric approach, based on linear transformations, can be used in solving initial and boundary value problems for linear differential equations as well.

**REFERENCES**

[1] M. Friedman, M. Ming, and A. Kandel. (1998). Fuzzy linear systems. *Fuzzy Sets and Systems*, 96:201-209.
[2] M. Ma, M. Friedman, and A. Kandel. (2000). Duality in fuzzy linear systems. *Fuzzy Sets and Systems*, 109:55-58.
[3] X. Wang, Z. Zhong, and M. Ha. (2001). Iteration algorithms for solving a system of fuzzy linear equations. *Fuzzy Sets and Systems*, 119:121-128.
[4] R. Ezzati (2011). Solving fuzzy linear systems. *Soft computing*, 15:193-197.
[5] J.F. Yin and K. Wang (2009). Splitting iterative methods for fuzzy system of linear equations. *Computational Mathematics and Modeling*, 20:326-335.
[6] T. Allahviranloo. (2004). Numerical methods for fuzzy system of linear equations. *Applied Mathematics and Computation*, 155:493-502.
[7] T. Allahviranloo. (2005). Successive over relaxation iterative method for fuzzy system of linear equations. *Applied Mathematics and Computation*, 162:189-196.
[8] T. Allahviranloo. (2005). The Adomian decomposition method for fuzzy system of linear equations. *Applied Mathematics and Computation*, 163:553-563.
[9] S. Abbasbandy, R. Ezzati, and A. Jafarian. (2006). LU decomposition method for solving fuzzy system of linear equations. *Applied Mathematics and Computation*, 172:633-643.
[10] T. Allahviranloo, E. Ahmady, N. Ahmady, and Kh. Shams Alketaby. (2006). Block Jacobi two stage method with Gauss-Sidel inner iterations for fuzzy systems of linear equations. *Applied Mathematics and Computation*, 175:1217-1228.
[11] M. Dehghan and B. Hashemi. (2006). Iterative solution of fuzzy linear systems. *Applied Mathematics and Computation*, 175:645-674.
[12] K. Wang and B. Zheng. (2006). Inconsistent fuzzy linear systems. *Applied Mathematics and Computations*, 181: 973-981.
[13] B. Asady, S. Abasbandy, and M. Alavi. (2005). Fuzzy general linear systems. *Applied Mathematics and Computation*, 169:34-40.
[14] T. Allahviranloo and M.A. Kermani. (2006). Solution of a fuzzy system of linear equations. *Applied Mathematics and Computation*, 175:519-531.
[15] B. Zheng and K. Wang. (2006). General fuzzy linear systems. *Applied Mathematics and Computation*, 181:1276-1286.
[16] X. Sun and S. Guo (2009). Linear formed general fuzzy linear systems. *Systems Engineering – Theory & Practice,* 29:92-98.
[17] T. Allahviranloo, N. Mikaeilvand, and M. Barkhordary. (2009). Fuzzy linear matrix equation. *Fuzzy Optim. Decis. Making*, 8:165-177.
[18] C.-X. Wu and M. Ma. (1991). Embedding problem of fuzzy number space: Part I. *Fuzzy Sets and Systems*, 44:33-38.
[19] C.-X. Wu and M. Ma. (1992). Embedding problem of fuzzy number space: Part III, *Fuzzy Sets and Systems*, 46:281-286.
[20] N.A. Gasilov, Sh.G. Amrahov, A.G. Fatullayev, H.I. Karakaş, and Ö. Akın. (2009). Existence theorem for fuzzy number solutions of fuzzy linear systems. Proceedings. pp. 436-439. 1st International Fuzzy Systems Symposium (FUZZYSS'09). TOBB University of Economics and Technology. October 1-2, 2009. Ankara, Turkey. Edited by T. Dereli, A. Baykasoglu and I.B. Turksen.
[21] N. Gasilov, Ş.E. Amrahov, A.G. Fatullayev, H.I. Karakaş, and Ö. Akın. (2009). A geometric approach to solve fuzzy linear systems. http://arxiv.org/ftp/arxiv/papers/0910/0910.4049.pdf
[22] R.B. Kearfott. (1996). *Rigorous Global Search: Continuous Problems*. Kluwer Academic Publishers, The Netherlands.